\documentclass[11pt]{article}
\usepackage[utf8]{inputenc}
\usepackage{graphicx}
\usepackage{amsmath, amsthm, amssymb} 
\usepackage{hyperref} 
\usepackage{multicol}
\usepackage{multirow}
\usepackage{caption}
\usepackage{xcolor}
\usepackage{subcaption}
\usepackage[labelfont=bf]{caption}
\usepackage{pdflscape}
\usepackage[noend]{algpseudocode}

\usepackage{algorithm}

\newcommand{\myalgsheader}[0]

\algnewcommand{\IIf}[1]{\State\algorithmicif\ #1\ \algorithmicthen}
\algnewcommand{\EndIIf}{\unskip\ \algorithmicend\ \algorithmicif}
\algnewcommand{\IElse}[1]{\State\algorithmicelse\ #1\ \algorithmicthen}
\algnewcommand{\IfThenElse}[3]{
	\State \algorithmicif\ #1\ \algorithmicthen\ #2\ \algorithmicelse\ #3}

\newcommand{\Atld}{\widetilde{A}}
\newcommand{\Etld}{\widetilde{E}}

\newcommand{\Rhat}{\widehat{R}}

\newcommand{\lvvert}{\hspace{1pt}\vert}
\newcommand{\rvvert}{\vert\hspace{1pt}}

\begin{document}
	\title{Mixed Precision FGMRES-Based Iterative Refinement for Weighted Least Squares}
	
	\author{Erin Carson\footnote{Department of Numerical Mathematics, Faculty of Mathematics and Physics, Charles University, \{carson, oktay\}@karlin.mff.cuni.cz. Both authors are supported by the Charles University GAUK project No. 202722, the Exascale Computing Project (17-SC-20-SC), a collaborative effort of the U.S. Department of Energy Office of Science and the National Nuclear Security Administration, and by the European Union (ERC, inEXASCALE, 101075632). Views and opinions expressed are those of the
authors only and do not necessarily reflect those of the European
Union or the European Research Council. Neither the European Union nor
the granting authority can be held responsible for them.} \hspace{1pt} and Eda Oktay\footnotemark[1]}
	\date{}
	\maketitle

	\paragraph{Abstract}
	With the recent emergence of mixed precision hardware, there has been a renewed interest in its use for solving numerical linear algebra problems fast and accurately. The solution of least squares (LS) problems $\min_x\|b-Ax\|_2$, where $A \in \mathbb{R}^{m\times n}$, arise in numerous application areas. Overdetermined standard least squares problems can be solved by using mixed precision within the iterative refinement method of Bj\"{o}rck, which transforms the least squares problem into an $(m+n)\times(m+n)$ ''augmented'' system. It has recently been shown that mixed precision GMRES-based iterative refinement can also be used, in an approach termed GMRES-LSIR. In practice, we often encounter types of least squares problems beyond standard least squares, including weighted least squares (WLS), $\min_x\|D^{1/2}(b-Ax)\|_2$, where $D^{1/2}$ is a diagonal matrix of weights. In this paper, we discuss a mixed precision FGMRES-WLSIR algorithm for solving WLS problems using two different preconditioners. 
	
	\section{Introduction}\label{sec:intro}
 Consider the weighted least squares problem
\begin{equation}
\min_x \Vert D^{1/2}(Ax-b) \Vert_2,
\label{eq:wls}
\end{equation}
where $A$ is an $m\times n$ matrix with $m\geq n$ and $D^{1/2}$ is a diagonal matrix of weights. In applications where the weights vary significantly in magnitude, transformation of the above into a standard least squares problem is not numerically stable \cite[Sec. 4.1.1]{bjorck1996numerical}. Note that in some applications, the weight matrix may be very ill-conditioned, for example, electrical networks, certain finite element problems, and interior point methods \cite[Remark 4.4.1]{bjorck1996numerical}.

Mathematically, the solution of \eqref{eq:wls} is given by the normal equations
\[
A^T D A x = A^T D b. 
\]
Letting $y=D(b-Ax)$, this is also equivalent to the augmented system
\begin{equation}
\begin{bmatrix}
\alpha D^{-1} & A \\ A^T & 0
\end{bmatrix}
\begin{bmatrix}
    \alpha^{-1} y \\ x
\end{bmatrix}
=
\begin{bmatrix}
    b \\ 0
\end{bmatrix},
\label{eq:augwls}
\end{equation}
where the scaling factor $\alpha$ has been introduced for stability; see \cite[Sec. 4.3.2 and Sec. 2.5.3]{bjorck1996numerical}. 

One potential way to solve weighted problems is based on Householder QR factorization, although one must use a Householder QR algorithm with row and column pivoting in order to guarantee a good accuracy \cite[Sec. 4.3.3]{bjorck1996numerical}. Powell and Reid investigated this in 1969 \cite{powell1969applying}, and a reworked backward error analysis in modern notation is given by Cox and Higham \cite{cohi98}. 

Another approach to solving such problems is using a mixed precision iterative refinement (IR) algorithm. With the advent of modern GPUs which feature multiple different hardware precisions, this approach has seen renewed interest. See Table \ref{tab:eps} for various IEEE precisions and their units roundoff. The general scheme of IR introduced in \cite{w:63} is given in Algorithm \ref{alg:ir}. Depending on the precision $u_f$ chosen for factorization $u_r$ for the residua computation, $u_s$ for the correction solve, and the working precision $u$ for storing data and the solution, one gets a variant of the mixed precision IR algorithm. To solve LS problems, the author in \cite{b:67} used a 2-precision IR approach to solve the augmented system; we call this LS-IR. To reduce the computation and memory cost, LS-IR never forms the augmented system explicitly, using only the QR factors of $A$. However, the QR factorization can be very expensive if $A$ is very large. Thus, the authors in \cite{carson2020three} use three precisions in LS-IR to further reduce the computation cost. Using a lower precision $u_f$ in the QR factorization can provide a cheaper algorithm while maintaining accuracy. Based on the analysis in \cite{carson2020three}, we can say as long as $\kappa_\infty(\Atld)\lesssim u_f^{-1}$, the backward error of three precision LS-IR is $\mathcal{O}(u)$ and the forward error is
\begin{equation}\label{eq:ferr_lsir}
    \dfrac{\|\Tilde{x}-\hat{\Tilde{x}}\|_\infty}{\|\Tilde{x}\|_\infty} \approx u_r cond(\Tilde{A},\Tilde{x})+u.
\end{equation}
Above, \eqref{eq:ferr_lsir} shows that if $u_f$ is chosen to be fp16, LS-IR converges only for $\kappa_\infty(\tilde{A})\lesssim 4.88\cdot 10^{4}$, which restricts the set of problems for which the algorithm is guaranteed to converge. 

\begin{algorithm}[htbp!]
	\caption{IR \cite{w:63} \label{alg:ir}} 
	\begin{algorithmic}[1]
		\Require{$A\in \mathbb{R}^{m \times n}$, $(m\geq n)$, $b\in \mathbb{R}^m$}
		\State{$Ax_0 = b$ \hfill in precision $u_f$; store in precision $u$} 
		\For{$i = 0$: $i_{max} - 1$}    
		\State{$r_i = b - Ax_i$ \hfill in precision $u_r$; store in precision $u$}
            \State{$Ad_{i+1} = r_i$ \hfill in precision $u_s$; store in precision $u$} \label{solve}
            \State{$x_{i+1} = x_i + d_{i+1}$ \hfill in precision $u$}
            \If{converged} {return $x_{i+1}$}
            \EndIf
		\EndFor
	\end{algorithmic}
\end{algorithm}

\begin{table}[]
	\centering
	\caption{Various IEEE precisions and their units roundoff. }
	\label{tab:eps}
	\begin{tabular}{cc}
		\hline
		\multicolumn{1}{c}{Precision} & \multicolumn{1}{c}{Unit Roundoff} \\ \hline
		fp16 (half) & $4.88\cdot 10^{-4}$ \\ 
		fp32 (single) & $5.96\cdot 10^{-8}$ \\ 
		fp64 (double) & $1.11\cdot 10^{-16}$ \\ 
		fp128 (quad) & $9.63\cdot 10^{-35}$ \\ \hline
	\end{tabular}
\end{table}

\section{GMRES-based approach}
To allow the use of low precision factorization for more ill-conditioned systems, we can adapt the GMRES-based approach of \cite{carson2020three} for non-weighted least squares problems, called GMRES-LSIR. The key difference with LS-IR is that this approach uses preconditioned GMRES to solve the linear system in line \ref{solve} of Algorithm \ref{alg:ir}.

For the non-weighted case, i.e., $D=I$, the authors in \cite{carson2020three} propose the left preconditioner 
\begin{equation}\label{eq:precond_lsir}
    M = \begin{bmatrix}
        \alpha I & Q\Rhat\\
        \Rhat^TQ^T& 0
    \end{bmatrix},
\end{equation}
composed of the QR factors of $A$ denoted as $Q$ and $\Rhat$.

According to \cite{carson2020three}, as long as $\kappa_\infty(A) \leq u^{-1/2}u_f^{-1}$, GMRES-LSIR has $\mathcal{O}(u)$ backward and forward error. In the simplified case where $u = u_r$, the forward error is \eqref{eq:ferr_lsir} if $\kappa_\infty(A)\leq u^{-1/3}u_f^{-2/3}$. Using the preconditioner $M$, the condition number of the preconditioned augmented system $M^{-1}\tilde{A}$ can be bounded by
\begin{equation}
    \kappa_\infty(M^{-1}\tilde{A}) \lesssim (1+2m\sqrt{n}\tilde{\gamma}_{mn}^f\kappa_\infty(A))^2, \quad \text{with } \quad \tilde{\gamma}_{mn}^f = \dfrac{cmn}{1-mnu_f}
\end{equation}
for a small constant $c$. This bound shows that even if $\kappa_\infty(A) \gg u_f^{-1}$, $M$ computed in precision $u_f$ will reduce $\kappa_\infty(M^{-1}\tilde{A})$, which is one of the most crucial benefits of GMRES-LSIR, allowing more ill-conditioned problems to be solved.

To solve WLSP using the iterative refinement-based approach, we can generally follow the same approach as described above. As in \cite[Sec. 3.1]{carson2020three}, we want to develop a preconditioner $M$ for our approach. Our aim is thus to find an effective and inexpensive preconditioner $M$ and prove that the resulting preconditioned coefficient matrix is sufficiently well-conditioned to guarantee the convergence of iterative refinement. In this study, we focus on two different preconditioners: a left preconditioner, which is a direct extension of the above, and a split preconditioner. 

\section{Left QR Preconditioner}
Consider the preconditioner
\begin{equation}
M_l=\begin{bmatrix}
    \alpha D^{-1} & Q\Rhat \\ \Rhat^TQ^T & 0
\end{bmatrix},
\label{eq:M}
\end{equation}
which is the direct extension of the $M$ in \eqref{eq:precond_lsir}, except with a $D^{-1}$ instead of $I$ in the $(1,1)$ block. We can explicitly write down the inverse of this matrix as 
\begin{equation*}
    M_l^{-1} = \begin{bmatrix}
    \frac{1}{\alpha}\left(D - DQ(Q^T D Q)^{-1} Q^T D \right) & DQ(Q^TDQ)^{-1}\Rhat^{-T} \\ \Rhat^{-1}(Q^TDQ)^{-1}Q^TD & -\alpha \Rhat^{-1}(Q^TDQ)^{-1}\Rhat^{-T}
\end{bmatrix}.
\end{equation*}

Let
\begin{equation*}
    \Etld = \Atld - M_l = \begin{bmatrix}
    0 & -E \\ -E^T & 0
\end{bmatrix},
\end{equation*}
where the error term $E$ is defined as 
\begin{equation}
    A+E = Q\Rhat,
    \label{eq:qrwls}
\end{equation}
due to finite precision computation of the QR factorization in some precision $u_f$ \cite{cohi98}.

Now, note that we can write 
\begin{align*}
M_l^{-1}\Atld &= M_l^{-1}(M_l+\Etld) = I+M_l^{-1}\Etld\\
\Atld^{-1}M_l &= (M_l+\Etld)^{-1}M_l \approx I-M_l^{-1}\Etld.
\end{align*}
Thus 
\begin{equation}
\kappa_\infty(M_l^{-1}\Atld) = \|M_l^{-1}\Atld\|_\infty \| \Atld^{-1}M_l\|_\infty \lesssim (1+\|M_l^{-1}\Etld\|_\infty)^2,
\label{eq:MAt}
\end{equation}
so if we bound the quantity $\|M_l^{-1}\Etld\|_\infty$, we obtain a bound on the condition number of the preconditioned system. 
We can explicitly write  
\begin{equation*}
   M_l^{-1}\Etld = \begin{bmatrix}
    -DQ(Q^TDQ)^{-1} \Rhat^{-T}E^T & -\frac{1}{\alpha}\left(D-DQ(Q^TDQ)^{-1}Q^TD \right)E \\ \alpha \Rhat^{-1}(Q^TDQ)^{-1}\Rhat^{-T} E^T & -\Rhat^{-1}(Q^TDQ)^{-1}Q^TDE
\end{bmatrix}. 
\end{equation*}
Note that if $D=I$ this reduces to the case in \cite{carson2020three} and gives \eqref{eq:precond_lsir}. 

The quantity $(QDQ)^{-1}Q^TD$ is what is called the ``scaled'' or ``weighted'' pseudoinverse \cite{stewart1989scaled}; see also \cite{o1990bounds}, \cite{wei1995upper}. In essence, the results of Stewart \cite{stewart1989scaled} show that this quantity is bounded and is independent of $D$. We can, however, still give a concrete bound on this quantity, which may in practice be a significant overestimate. Let $Q_D = D^{1/2}Q$. Then 
\begin{eqnarray}
\Vert (Q^T D Q)^{-1} Q^T D \Vert_2 &= \Vert (Q_D^T Q_D)^{-1} Q_D^T D^{1/2} \Vert_2 \nonumber \\
&\leq \Vert Q_D^\dagger \Vert_2 \Vert D^{1/2} \Vert_2 \nonumber \\
&\leq \Vert Q^\dagger \Vert_2 \Vert D^{-1/2} \Vert_2 \Vert D^{1/2} \Vert_2 \nonumber\\
&= \kappa_2^{1/2}(D). \label{eq:spinvbound}
\end{eqnarray}
As an alternative, we can let 
\begin{equation*}
    \Vert (Q^T D Q)^{-1} Q^T D \Vert_2 \equiv \rho,
\end{equation*}
and then argue that this $\rho$ is hopefully of reasonable size since it is independent of $D$, using the results of Stewart (we will do this for now). 

Now, we want to derive a bound for $\| M_l^{-1}\Etld\|_\infty$. Following \cite{carson2020three}, we have 
\begin{align}
\| M_l^{-1} \Etld \|_\infty &\leq \max (\| (M_l^{-1}\Etld)(1,1)\|_\infty +   \| (M_l^{-1}\Etld)(1,2)\|_\infty, \nonumber \\
&\phantom{\leq}\qquad\qquad\qquad, \| (M_l^{-1}\Etld)(2,1)\|_\infty +   \| (M_l^{-1}\Etld)(2,2)\|_\infty )\nonumber \\
&\leq \sqrt{m} \max \left( \| (M_l^{-1}\Etld)(1,1)\|_F, \| (M_l^{-1}\Etld)(2,1)\|_F \right) \nonumber \\
&\phantom{\leq}+ \sqrt{n} \max \left( \| (M_l^{-1}\Etld)(1,2)\|_F, \| (M_l^{-1}\Etld)(2,2)\|_F \right),\label{eq:MinvE}
\end{align}
where $(M_l^{-1}\Etld)(i,j)$ denotes the $(i,j)$-block of $M_l^{-1}\Etld$. 

Using $|\alpha|\approx \Vert A^{\dagger} \Vert_2^{-1}$ and ignoring terms of order $u_f^2$, 
\begin{align*}
\| (M_l^{-1}\Etld)(1,1)\|_F &= \| DQ(Q^TDQ)^{-1} \Rhat^{-T}E^T\|_F \leq \rho \Vert A^\dagger \Vert_2 \Vert E^T \Vert_F \\
\| (M_l^{-1}\Etld)(1,2)\|_F &= \left\Vert \frac{1}{\alpha}D \left(I-Q(Q^TDQ)^{-1}Q^TD \right)E\right\Vert_F \leq \rho \Vert D \Vert_2 \Vert A^\dagger \Vert_2 \Vert E \Vert_F  \\
\| (M_l^{-1}\Etld)(2,2)\|_F&= \| \Rhat^{-1}(Q^TDQ)^{-1}Q^TDE\|_F \leq \rho \Vert A^\dagger \Vert_2 \Vert E \Vert_F.
\end{align*}

For the final term, we have
\begin{align*}
\| (M_l^{-1}\Etld)(2,1)\|_F&= \| \alpha \Rhat^{-1}(Q^TDQ)^{-1}\Rhat^{-T} E^T\|_F \\
&= \|\alpha \Rhat^{-1} (Q^TDQ)^{-1}Q^TD(Q^TD)^{\dagger} \Rhat^{-T} E^T \|_F\\
&\leq \rho \Vert D^{-1}\Vert_2 \Vert A^\dagger \Vert_2 \Vert E^T\Vert_F. 
\end{align*}

Then plugging into \eqref{eq:MinvE}, and letting $\theta \equiv \max (1, \| D\|_2, \|D^{-1}\|_2)$, we have
\begin{eqnarray*}
\Vert M_l^{-1}\Etld\Vert_\infty & \leq (\sqrt{m}+\sqrt{n}) \rho \theta \Vert A^\dagger \Vert_2 \Vert E \Vert_F .
\end{eqnarray*}

Assuming that we have a standard QR factorization, $\Vert E \Vert_F \leq \tilde{\gamma}_{mn}^f \Vert A \Vert_F$, and thus
\begin{eqnarray}\label{eq:err_wlsir_qr}
\Vert M_l^{-1}\Etld\Vert_\infty & \leq (\sqrt{m}+\sqrt{n}) \rho \theta \Vert A^\dagger \Vert_2 \Vert E \Vert_F\\
&\leq 2m\sqrt{n} \rho \theta \tilde\gamma_{mn}^f \kappa_\infty(A). 
\end{eqnarray}
Note that in the case that we use a standard QR factorization and $D=I$ (and thus $\rho=1$), we recover the same bound as in \cite[Section 3.1]{carson2020three}. 
Finally, plugging into \eqref{eq:MAt}, 
we obtain
\begin{equation}
\kappa_\infty(M_l^{-1}\Atld) \lesssim \left(1+2m\sqrt{n}\rho \theta \tilde\gamma_{mn}^f \kappa_\infty(A) \right)^2.
\label{eq:kinfbound}
\end{equation}

\section{Block-diagonal Split Preconditioner}
Note that the dependence of the bound in \eqref{eq:err_wlsir_qr} on $D$ is not ideal since $D$ can be very ill-conditioned. For the behavior of the condition number of $M_l^{-1}\Atld$, see Section \ref{sec:exp_ml}. To construct a preconditioner which results in a preconditioned matrix whose conditioning does not depend on $D$, we thus consider using the block diagonal preconditioner 
\begin{equation}\label{block}
    M_b = \begin{bmatrix}
        \alpha D^{-1} & 0\\
        0 & \Hat{C}
    \end{bmatrix},
\end{equation}
where $\Hat{C}$ is a symmetric positive definite approximation to the Schur complement, $\alpha^{-1}A^TDA$ \cite{rozlovznik2018saddle}.

Although $M_b$ can be used as a left/right preconditioner, to obtain a symmetric preconditioned system, we consider using it as a split preconditioner and study the system
\begin{align*}
    M_b^{-1/2}\Tilde{A}M_b^{-1/2} &= \begin{bmatrix}
       (\alpha D^{-1})^{-1/2} & 0\\
        0 & \Hat{C}^{-1/2} 
    \end{bmatrix}\begin{bmatrix}
       \alpha D^{-1} & A\\
        A^T & 0 
    \end{bmatrix}\begin{bmatrix}
       (\alpha D^{-1})^{-1/2} & 0\\
        0 & \Hat{C}^{-1/2} 
    \end{bmatrix} \\&= \begin{bmatrix}
       (\alpha D^{-1})^{-1/2}(\alpha D^{-1})(\alpha D^{-1})^{-1/2} & (\alpha D^{-1})^{-1/2}A^T\Hat{C}^{-1/2}\\
        \Hat{C}^{-1/2}A^T(\alpha D^{-1})^{-1/2}& 0
    \end{bmatrix} \\
    &= \begin{bmatrix}
        I & \Hat{A}\\
        \Hat{A}^T & 0
    \end{bmatrix}.
\end{align*}

No applicable analysis of the forward and backward errors in split preconditioned GMRES exists, but there is an existing analysis of these quantites for a general split preconditioned flexible GMRES method (FGMRES). Thus, to be able to use split preconditioners and discuss their effects on the stability of this general approach, we use FGMRES instead of GMRES in GMRES-LSIR. Our new variant to solve WLSP is  therefore called FGMRES-WLSIR.

We use a symmetric positive definite approximation $\Hat{C}$ such that it is spectrally equivalent to the matrix $\alpha^{-1}A^TDA$, i.e., there exist positive constants $0<\hat{\gamma}\leq \hat{\delta}$ such that 
\begin{equation*}
    \hat{\gamma}(\Hat{C}y,y)\leq (B^T(\alpha D^{-1})^{-1}Ay,y)\leq \hat{\delta}(\Hat{C}y,y)
\end{equation*}
for all vectors $y\in \mathbb{R}^n$.
Using \cite[Proposition 3.4]{rozlovznik2018saddle}, we can find the spectrum of the preconditioned matrix:
\begin{equation*}
    sp(M_b^{-1/2}\Tilde{A}M_b^{-1/2}) = \frac{1}{2}\left(1\pm \sqrt{1+4\sigma^2_k(\hat{A})}\right) \quad for \quad k=1,\ldots,n-r,
\end{equation*}
where rank($\hat{A}$) = $n-r$ and $0<\sigma_{n-r}(\hat{A})\leq \cdots\leq \sigma_1(\hat{A})$ are the singular values of $\hat{A}$. 

We assume that the Schur complement is computed exactly and we consider the left-preconditioned matrix
\begin{equation*}
    M_b^{-1}\Tilde{A} = \begin{bmatrix}
        I&\alpha^{-1}A^TD\\
        \alpha^{-1}A(A^TDA)^{-1}&0
    \end{bmatrix},
\end{equation*}
which is a nonsymmetric diagonalizable matrix with three distinct eigenvalues $\left\{1,\frac{1}{2}(1\pm \sqrt{5})\right\}$. This makes the block diagonal preconditioner a suitable choice for FGMRES since the method can converge in a small number of iterations. However, our experiments show that $\kappa_\infty(M_b^{-1}\Tilde{A})$ can be very large, often larger than $\kappa_\infty(\Tilde{A})$, when $A$ is ill-conditioned. An ill-conditioned preconditioned matrix makes the preconditioner unsuitable for proving the backward stability of FGMRES although we show in Section \ref{sec:exps} that in practice, the split preconditioner improves the condition number of the preconditioned matrix.

Using the analysis in \cite{rozlovznik2018saddle}, we can obtain the bound 
\begin{equation}\label{eq:cond_block}
		\kappa_\infty(M_{b}^{-1/2}\tilde{A}M_{b}^{-1/2})\leq \dfrac{\lvvert1+\sqrt{1+4\sigma_1^2(\hat{A})}\rvvert}{\lvvert1-\sqrt{1+4\sigma_{n-r}^2(\hat{A})}\rvvert}(n+m),
	\end{equation}
where $0<\sigma_{n-r}(\hat{A})\leq \ldots \leq \sigma_1(\hat{A}),$
\begin{equation*}
	\hat{A} = \dfrac{1}{\sqrt{\alpha}}D^{1/2}A^T\hat{C}^{-1/2} \text{ with } \hat{C} \approx \alpha^{-1}A^TDA.
\end{equation*}

Although there is still a dependence on $D$ in \eqref{eq:cond_block}, we can numerically eliminate its impact by computing the Schur complement in a specific manner. For $\hat{C}$, we compute $D^{1/2}A$ in high precision and use the R-factor (from lower precision QR factorization) of it in $M_b^{1/2}$. The numerical experiments in Section \ref{sec:exps} show that this trick reduced the effect of the conditioning of $D$.

On the other hand, using a split preconditioner in (F)GMRES has several disadvantages regarding error analysis. Using the analysis in \cite{cd:23}, we see that for the block diagonal split preconditioner (the left and right preconditioners are both $M_b^{1/2}$), the forward error of FGMRES is bounded by
\begin{equation}
	\dfrac{\|x-\bar{x}_k\|}{\|x\|}\lesssim \kappa_\infty(M_{b}^{-1/2}\tilde{A}M_{b}^{-1/2})\kappa_\infty(M_{b}^{1/2})\mathcal{O}(u),	
 \label{eq:ferrfgmres}
\end{equation}
whereas for any left preconditioner, the bound is given as
\begin{equation*}
	\dfrac{\|x-\bar{x}_k\|}{\|x\|}\lesssim \kappa_\infty(M_{l}^{-1}\tilde{A})\mathcal{O}(u).	
\end{equation*}
Even though $M_{b}^{-1/2}\tilde{A}M_{b}^{-1/2}$ is well-conditioned, $M_{b}^{1/2}$ can still be ill-conditioned, which can cause the forward error to be more than $1$ and thus cannot guarantee FGMRES-WLSIR convergence. This result shows the need for a bound on $\kappa_\infty(M_{b}^{1/2})$.

\section{Numerical Experiments}\label{sec:exps}
To illustrate the impact of different preconditioners on FGMRES-WLSIR convergence, we perform several numerical experiments in MATLAB. The experiments are performed on a computer with AMD Ryzen 5 4500U having 6 CPUs and 8 GB RAM with OS system Ubuntu 22.04 LTS. In our numerical experiments, we used built-in MATLAB datatypes for double and single precisions. To simulate half-precision floating point arithmetic, we use the chop library and associated functions from \cite{hp:19} available at \verb|https://github.com/higham/chop| and \verb|https://github.com/SrikaraPranesh/|\\ \verb|LowPrecision_Simulation|. Code for our FGMRES-WLSIR 
and associated functions can be found in the repository \verb|https://github.com/edoktay/fgmreswlsir|.

For the examples in this section, we set $A' =$ \texttt{gallery(`randsvd', [100,10], 1e2, 3)}, which creates a matrix with 2-norm condition number $10^2$ of size $100\times 10$ with geometrically distributed singular values. We prescale the rows of $A'$ so that they have drastically different sizes. We use 9 different scalings $A = S^{-1}A'$ where $S$ is a diagonal matrix created by \texttt{diag(logspace(1, j, 100)} where $j\in[1,2,4,6,8,10,12,14,16]$. We set the weight matrix $D$ such that each row $i$ of $A$ has $\max |A(i,j)|=1$. The scaling factor $\alpha = 2^{-1/2} \sigma_n$, where $\sigma_n$ is the smallest singular value of $A$. We compute QR factorizations in half, single, and double precision and construct the corresponding preconditioners $M$. We then measure the infinity-norm condition number of the preconditioned system. 

In each figure, the condition number of the preconditioned systems is represented as colored lines. Each color shows half (red), single (green), and double (blue) precisions used for computing QR factorizations for the preconditioners. The dashed black line gives the condition number of the unpreconditioned augmented system and the dotted black line gives the inverse of the unit roundoff for the FGMRES-WLSIR working precision $u$. The figure shows that the convergence of FGMRES-WLSIR is guaranteed only when the colored solid lines remain below the dashed line. Numerically, this shows the cases when $\kappa_\infty(M^{-1}\Atld)\leq u^{-1}$. Only in this case can we guarantee that the forward error of FGMRES is less than $1$ and thus FGMRES-WLSIR converges.

\subsection{Dependence of $\Vert M_l^{-1}\Etld\Vert_\infty$ on $D$}\label{sec:exp_ml}
\begin{figure}[htb!]
\centering
\includegraphics[trim={3cm 8cm 4cm 8cm}, clip, width=.49\textwidth]{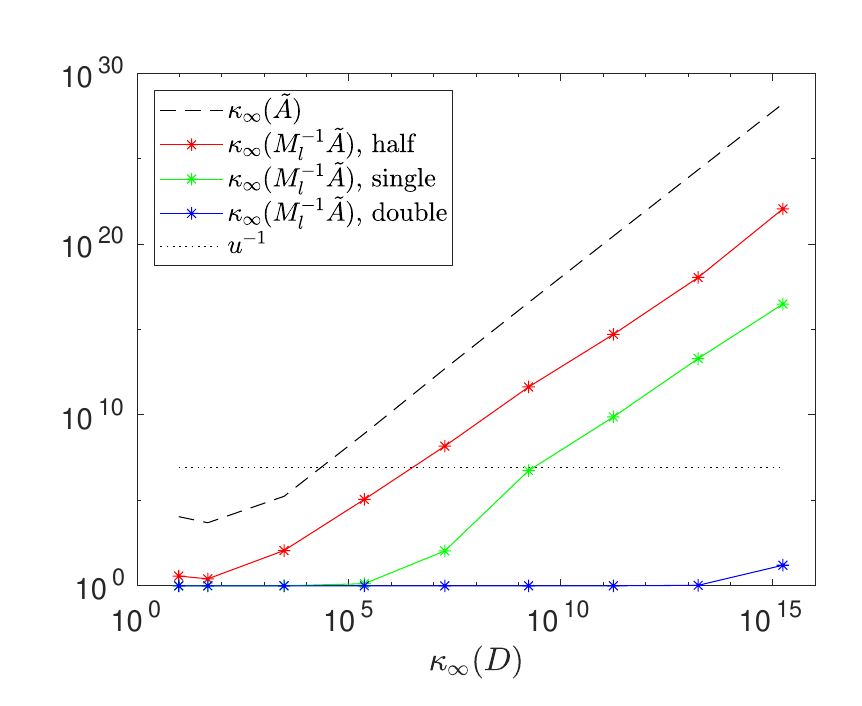}
\caption{Measured condition number of the preconditioned systems and estimates of the bound on the condition number \eqref{eq:kinfbound} where the preconditioners \eqref{eq:M} are constructing using QR factorizations computed in various precisions, versus the condition number of the weight matrix $D$. The working precision for FGMRES-WLSIR is assumed to be single precision.}
\label{fig:condbound}    
\end{figure}

We perform a numerical experiment to demonstrate how $\kappa_\infty(M_l^{-1}\Atld)$ changes with the conditioning of $D$. For this example, we assume that FGMRES-WLSIR is performed using single precision as the working precision $u$. The results are shown in Figure \ref{fig:condbound}.

The figure shows that using $M_l$ in double precision preserves stability even if $D$ is very badly-conditioned. When we switch to single precision, we see $M_l$ is useful (convergence of FGMRES is guaranteed) when $\kappa_\infty(D)< 10^{10}$. On the other hand, using half-precision limits the usability of the preconditioner since the condition number of $D$ should be at most $10^6$, which may not be the case in real-world problems. We thus conclude that using lower precision can be tricky for $M_l$ due to the fact that its error bound depends on the conditioning of $D$. We also observe from the black dashed line that half-precision $M_l$ does not change the conditioning of the preconditioned system significantly.

\subsection{$M_l$ versus $M_b$}\label{sec:exp_mb}
To examine the effect of preconditioners on the conditioning of the augmented system, we construct Table \ref{tab:cond} using random dense matrices with a randomly generated solution vector $b$. We construct matrix $A$ in a similar manner as in Section \ref{sec:exp_ml}.
The table shows that the preconditioner $M_l$ does not decrease the condition number sufficiently, whereas $M_b$ works well. We finally observe from the last column that even though $M_b^{1/2}$ is ill-conditioned since the preconditioned system is very well conditioned, FGMRES-WLSIR still converges due to the forward error constraint in \eqref{eq:ferrfgmres}. 

\begin{table}[]
\centering
	\caption{Condition numbers of $\Atld$, right preconditioner, and preconditioned augmented matrices.}
	\label{tab:cond}
	\begin{tabular}{ccccc}
		\hline
		$\kappa_2(A)$ &$\kappa_\infty(\Atld)$ &$\kappa_\infty(M_l^{-1}\Atld)$ & $\kappa_\infty(M_b^{-1/2}\Atld M_b^{-1/2})$& $\kappa_\infty(M_b^{1/2})$ \\ \hline
                1.00e+02& 1.12e+04 & 3.73e+00 & 8.56e+01 & 1.05e+03\\
                1.12e+02& 4.92e+03 & 2.55e+00 & 7.80e+01 & 3.96e+02\\
		      1.47e+02	&1.73e+05 & 1.16e+02 & 5.79e+01& 3.60e+02 \\ 
			1.91e+02	& 8.04e+08 & 1.14e+05 & 4.03e+01& 3.08e+04\\
                2.47e+02	& 5.06e+12 & 1.50e+08 & 3.02e+01& 2.48e+06\\ 
			3.21e+02	& 3.81e+16 & 4.36e+11 & 2.62e+01& 2.16e+08\\
			4.22e+02	&2.97e+20 &5.24e+14 &2.30e+01& 1.87e+10 \\ 
			5.61e+02	&2.26e+24 &1.12e+18 &2.05e+01& 1.63e+12 \\ 
			7.59e+02	&1.67e+28 &1.16e+22 &1.86e+01& 1.36e+14\\ \hline
	\end{tabular}
\end{table}

To demonstrate the difference in the numerical behavior of the preconditioned system with $M_l$ and $M_b$ with the conditioning of $D$, we also perform several numerical experiments using \verb|ash958| and \verb|robot24c1_mat5| matrices from the SuiteSparse collection \cite{suitesparse}. We prescale the rows of $A'$ so that they have drastically different sizes. We use 9 different scalings $A = SA'$ where $S$ is the scaling matrix used in Section \ref{sec:exp_ml}. The properties of both matrices are given in Table \ref{tab:matrices}. For the experiments in this section, we assume that we use double precision as the working precision $u$ in FGMRES-WLSIR. 

\begin{table}[htb!]
    \centering
	\caption{Properties of matrices from the SuiteSparse collection. }
	\label{tab:matrices}
	\begin{tabular}{ccccc}
		\hline
		Name & $m$ &$n$ &$\kappa_2(A)$& $\#$nnz \\ \hline
		\verb|ash958| & 958 & 292 & 3.2014&1916\\
\verb|robot24c1_mat5| & 404 &302 & $3.33\times 10^{11}$& 15118 \\ \hline
	\end{tabular}
\end{table}

\begin{figure}[htb!]
    \centering
    \includegraphics[trim={3cm 8cm 4cm 8cm}, clip, width=.49\textwidth]{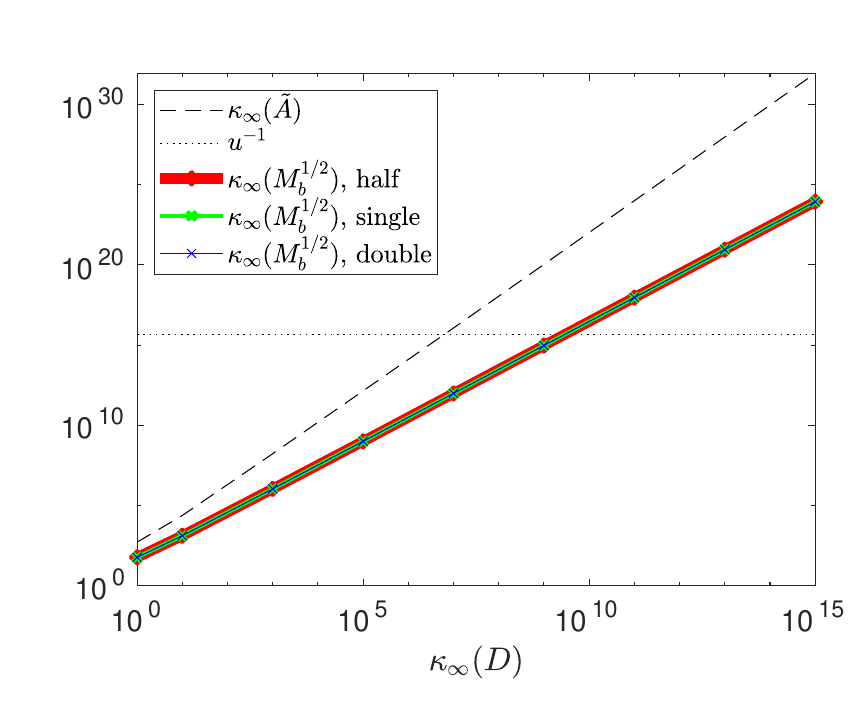}
    \includegraphics[trim={3cm 8cm 4cm 8cm}, clip, width=.49\textwidth]{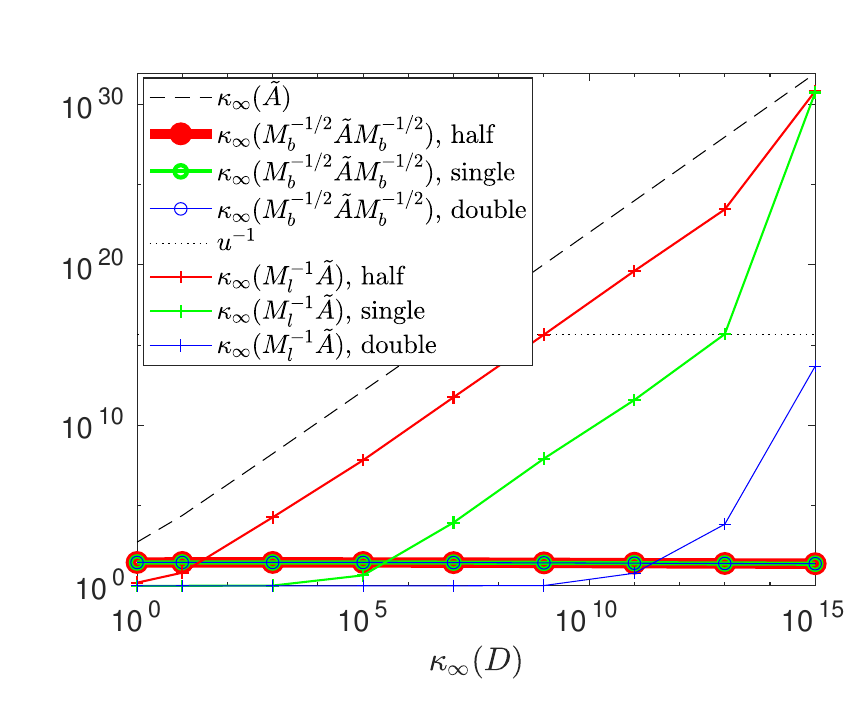}
    \caption{Measured condition number of the preconditioners (left) and preconditioned systems (right) using \textit{ash958} matrix as $A$ where $M_l$ and $M_b$ are constructed using QR factorizations in various precisions, versus the condition number of the weight matrix $D$. The working precision for FGMRES-WLSIR is assumed to be double precision.}
    \label{fig:ash958}
\end{figure}

We observe from the right plot in Figure \ref{fig:ash958} that single precision $M_l$ can handle up to $\kappa_\infty(D)< 10^{13}$ however, we note that the algorithm can still be expensive due to QR factorization. On the other hand, we see from the same plot that using $M_b$, even with a very ill-conditioned $D$ will provide a well-conditioned preconditioned coefficient matrix. Numerically, we don't see the effect of $D$ because of our way of computing the Schur complement. However, because of the dependence of the forward error of FGMRES on $\kappa_\infty(M_b^{1/2})$, we need to look at the left plot as well. From the left plot we see that even if numerically, the preconditioned system is very well-conditioned, because of the conditioning of the right split preconditioner (in our case, it is equivalent to the left split preconditioner), $\kappa_\infty(M_b^{1/2})$ is sufficiently small only when $\kappa_\infty(D)<10^9$. For \verb|ash958|, we can say that using the split block diagonal preconditioner $M_b$ does not give any significant advantage over $M_l$ in fp16.

\begin{figure}[htb!]
    \centering
    \includegraphics[trim={3cm 8cm 4cm 8cm}, clip, width=.49\textwidth]{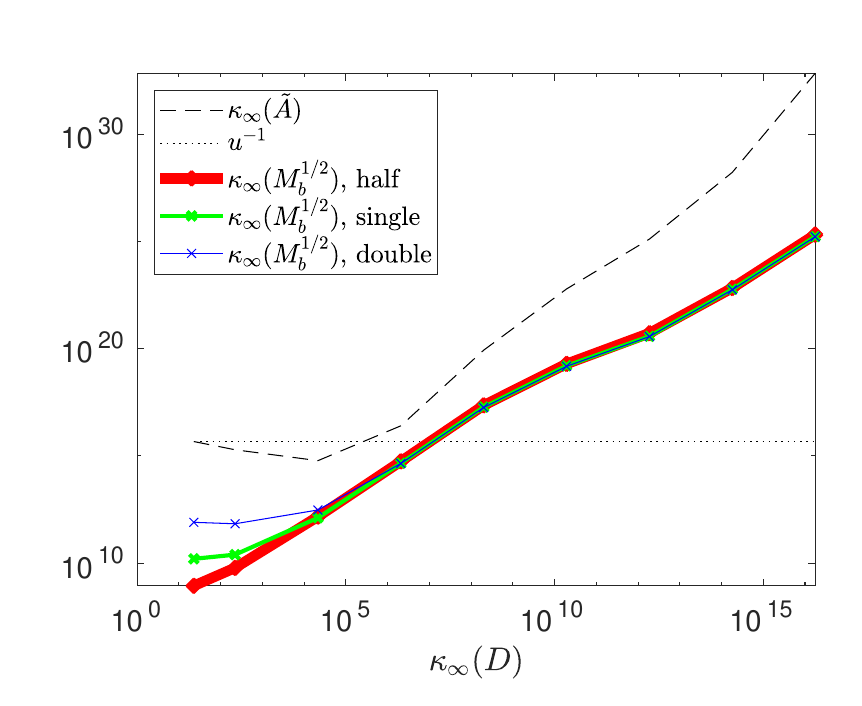}
    \includegraphics[trim={3cm 8cm 4cm 8cm}, clip, width=.49\textwidth]{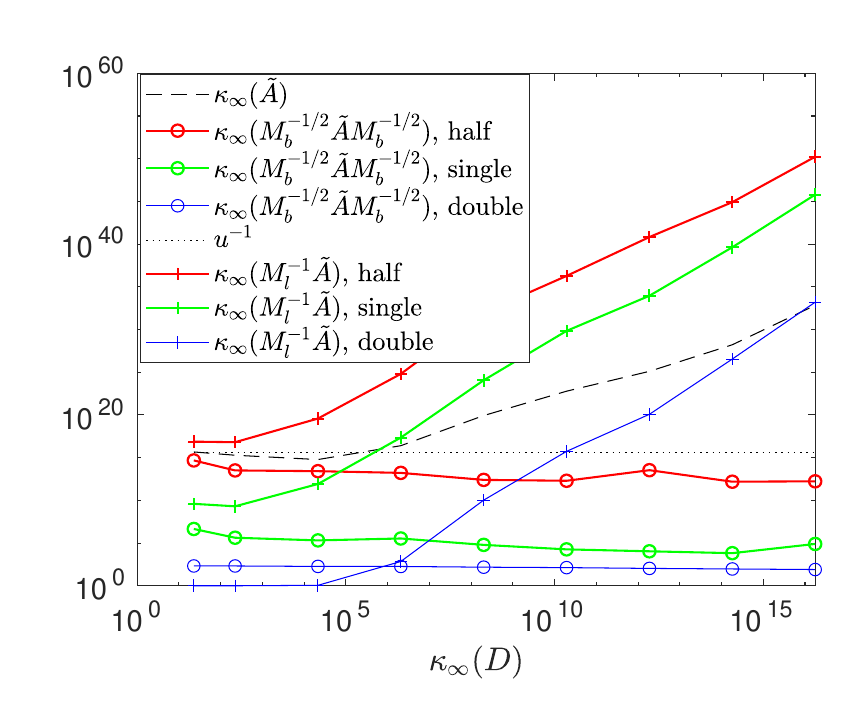}
    \caption{Measured condition number of the preconditioners (left) and preconditioned systems (right) using \textit{robot24c1\_mat5} matrix as $A$ where $M_l$ and $M_b$ are constructed using QR factorizations in various precisions, versus the condition number of the weight matrix $D$. The working precision for FGMRES-WLSIR is assumed to be double precision.}
    \label{fig:robot}
\end{figure}

From the right plot in Figure \ref{fig:robot}, we see that using half-precision in computing $M_l$ is not useful even if $D$ is very well-conditioned. Moreover, the left preconditioner works only when $\kappa_\infty(D)<10^7$. For the split preconditioner, we again observe the same behavior as the previous example. From the left preconditioner, we again expect a limitation coming from the conditioning of $M_b^{1/2}$.  We observe that in any precision, $\kappa_\infty(M_b^{1/2})\leq u^{-1}$ when $\kappa_\infty(D)<10^7$. However, even though $\kappa_\infty(M_b^{1/2})\leq u^{-1}$, since the forward error bound is obtained via the multiplication of both condition numbers, it will be close to 1 for $M_b$ in fp16 and fp32 when $D$ is well-conditioned. Therefore, using fp16 in computing both preconditioners is not applicable to this matrix. Furthermore, we see that for \verb|robot24c1_mat5|, $M_l$ is more useful than $M_b$ in terms of fp32 and fp64 applicability.

\section{Conclusion}
In various areas, problems may need a very badly-conditioned weight matrix to be able to be modeled as least squares problems. Most of the available fast least squares solvers are not directly usable for weighted cases. One thus needs to make some changes in algorithms such as changing the preconditioning. GMRES-LSIR is one of the potential iterative solvers in the literature for solving least squares problems fast and accurately using lower precision. With this motivation, our goal was to extend GMRES-LSIR to solve weighted least squares problems. For this extension, we examined the use of two different preconditioners, namely left and block split preconditioners, in the algorithm. To construct an error bound for the split preconditioner, we use FGMRES instead of the GMRES algorithm and introduce our approach, FGMRES-WLSIR.

Using the analyses in the literature, we introduce error bounds for both preconditioners under assumptions on the conditioning of $A$. From our analysis, we observe that the dependence of the weight matrix in the error bound of the left preconditioner limits the use of low precisions in our approach. To reduce this dependence numerically, we study the block split preconditioner. By computing the Schur complement in a specific way, we numerically reduce the dependence on the conditioning of $D$. However, from the preconditioned FGMRES analysis, we observe that the forward error of FGMRES also depends on the conditioning of the right preconditioner in the split preconditioner case. We therefore examine the conditioning of both preconditioners and their ranges of applicability numerically. We conclude that since the conditioning of the weight matrix is highly problem-dependent, we cannot generalize which preconditioner is more useful for FGMRES-WLSIR. Although our numerical experiments show that using the block split preconditioner in half-precision may be used in more ill-conditioned systems, in real applications, $D$ might be worse-conditioned. 

Because of the dependence of both preconditioners on $D$ and the dependence of split preconditioned FGMRES on the right preconditioner, further study is warranted. Future studies can focus either on the choice of a preconditioner or another iterative approach other than FGMRES. In any case, the optimal approach will ideally have error bounds independent of the weight matrix. 
	\bibliographystyle{siamplain}
	\bibliography{main_arxiv_new}

\end{document}